\documentclass[12pt,a4paper,oneside]{article}
\usepackage{amsmath,amsfonts,amssymb,latexsym}

\newtheorem{theorem}{Theorem}[section]
\newtheorem{lemma}[theorem]{Lemma}
\newtheorem{definition}[theorem]{Definition}
\newtheorem{proposition}[theorem]{Proposition}
\newtheorem{corollary}[theorem]{Corollary}
\newtheorem{remark}[theorem]{Remark}
\newtheorem{example}[theorem]{Example}
\newenvironment{proof}
{\bigskip\noindent{\sc Proof.}\ \ \rm }{\hfill$\Box$\bigskip}
\title{ On a new probabilistic representation for the solution of the heat equation }
\author{ Paolo Da Pelo \quad\quad Alberto Lanconelli }
\date{\empty}
\begin{document}
\maketitle { \noindent \small

\!\!\!\!\!\!\!\!\!\!\!\!\!\!\!\!
\begin{center}
{ \noindent
\begin{tabular}{cc}
& Dipartimento di Matematica\\
& Universita' degli Studi di Bari\\
& Via E. Orabona, 4- CAP 70125 Bari, Italia\\
& E-mails: dapelo@dm.uniba.it, lanconelli@dm.uniba.it \\

\end{tabular}
}
\end{center}

} \numberwithin{equation}{section}

\bigskip

\begin{abstract}
We obtain a new probabilistic representation for the solution of the
heat equation in terms of a product for smooth random variables
which is introduced and studied in this paper. This multiplication,
expressed in terms of the Hida-Malliavin derivatives of the random
variables involved, exhibits many useful properties which are to
some extents opposite to some peculiar features of the Wick product.

\end{abstract}
     Key words and phrases: heat equation, Hida-Malliavin derivative, second quantization operators,
     smooth random variables.\\

AMS 2000 classification: 60H30, 60H07.

\section{Introduction}

In the theory of stochastic differential equations (SDEs), It\^o
SDEs driven by a one-dimensional Brownian motion are the most
investigated ones; they are of the form,
\begin{eqnarray}\label{1.1}
X_t^x(\omega)=x+\int_0^tb(X_s^x(\omega))ds+\int_0^t\sigma(X_s^x(\omega))dB_s(\omega),\quad
0\leq t\leq T,
\end{eqnarray}
where $b,\sigma:\mathbb{R}\to\mathbb{R}$ are measurable functions,
$\{B_t\}_{0\leq t\leq T}$ is a one-dimensional Brownian motion,
$x\in\mathbb{R}$ is the initial value and the last term in the
right-hand
side is an It\^o integral.\\
Stochastic differential equations are, in fact, integral equations;
the non-differentiability of the Brownian motion with respect to $t$
forces the integral interpretation of the equation. White noise
analysis is a theory of stochastic distributions, analogous in the
construction to the classical distribution theory; within this
framework one can give a precise meaning to the time derivative of
the Brownian motion and obtain the so called \emph{white noise}
$W_t(\omega)$. One of the crucial differences between ordinary and
stochastic differential equations is that the just mentioned
differentiation reduces equation (\ref{1.1}) to:
\begin{eqnarray*}
\frac{dX_t^x(\omega)}{dt}=b(X_t^x(\omega))+\sigma(X_t^x(\omega))\diamond
W_t(\omega).
\end{eqnarray*}
The symbol $\diamond$ denotes the \emph{Wick product} of
$\sigma(X_t^x(\omega))$ and $W_t(\omega)$. This operation, defined
for a certain class of stochastic distributions, is therefore the
hidden core of It\^o integration theory and It\^o
stochastic differential equations. \\

It is well known (see for instance \cite{HOUZ}) that for any $t\in
[0,T]$ and $n\geq 1$,
\begin{eqnarray*}
B_t^{\diamond n}&:=&\underbrace{B_t\diamond\cdot\cdot\cdot\diamond B_t}_{n-times}\\
&=&h_{n,t}(B_t),
\end{eqnarray*}
where $h_{n,t}$ denotes the $n$-th order Hermite polynomial with
parameter $t$ and leading coefficient one. More generally, if
$g:\mathbb{R}\to\mathbb{R}$ is a real analytic function with
expansion
\begin{eqnarray*}
g(x)=\sum_{n\geq 0}a_nx^n,\quad x\in\mathbb{R},
\end{eqnarray*}
then
\begin{eqnarray*}
u(t,x)&:=&g^{\diamond}(B_t)\vert_{B_t=x}\\
&=&\sum_{n\geq 0}a_nB_t^{\diamond n}\Big\vert_{B_t=x}
\end{eqnarray*}
solves the backward heat equation
\begin{eqnarray*}
\partial_{t}u(t,x)+\frac{1}{2}\partial_{xx}u(t,x)=0.
\end{eqnarray*}
See \cite{L1} and the references quoted there.\\

The aim of the present paper is to establish a similar
correspondence between a certain stochastic multiplication and the
solution of the forward heat equation
\begin{eqnarray}\label{1.2}
\partial_{t}u(t,x)-\frac{1}{2}\partial_{xx}u(t,x)=0.
\end{eqnarray}
More precisely we introduce a new product for smooth random
variables, defined in terms of Hida-Malliavin derivatives  and named
anti-Wick product, and we show that the unique solution of
(\ref{1.2}) with initial condition $f$ can be explicitly represented
in terms of this multiplication and the data $f$. We also show that
the solution of (\ref{1.2}) with initial data $f\cdot g$ can be
written as the anti-Wick product of the two solutions with data $f$
and
$g$, respectively.\\

The paper is organized as follows: Section 2 recalls basic
definitions, notations and theorems from the Malliavin calculus and
White Noise theory. In Section 3 we define the anti-Wick product for
smooth random variables and study some of its crucial properties, in
particular the representation stated in Theorem 3.7. In Section 4 we
present the main result of the paper, i.e. the new probabilistic
representation for the solution of the heat equation. Finally in the
Appendix we give two formulas that relate Wick and anti-Wick
products.

\section{Framework}
In this section we set the necessary tools for our investigation.
For more information we refer the reader to one of the books
\cite{HKPS},\cite{HOUZ},\cite{J},\cite{Kuo}, \cite{N},\cite{Ob}.\\
Let $(\Omega,\mathcal{F},\mathcal{P})$ be a complete probability
space which carries a one dimensional Brownian motion
$\{B_t\}_{0\leq t\leq T}$. Assume that $\mathcal{F}=\mathcal{F}_T$
where $\{\mathcal{F}_t\}_{0\leq t\leq T}$ denotes the augmented
Brownian filtration; as a consequence of this assumption we get that
the set
\begin{eqnarray*}
\Big\{\mathcal{E}(f):=\exp\Big\{\int_0^Tf(s)dB_s-\frac{1}{2}\int_0^Tf^2(s)ds\Big\},
f\in\mathcal{L}^2([0,T])\Big\}
\end{eqnarray*}
is total in $\mathcal{L}^2(\Omega,\mathcal{F},\mathcal{P})$
($\mathcal{L}^2(\Omega)$ for short).\\
According to the Wiener-It\^o chaos decomposition theorem any
$X\in\mathcal{L}^2(\Omega)$ can be uniquely represented as
\begin{eqnarray*}
X&=&\sum_{n\geq 0}I_n(h_n)\quad(\mbox{convergence in }
\mathcal{L}^2(\Omega)),
\end{eqnarray*}
where $I_0(h_0):=E[X]$ and for $n\geq 1$,
$h_n\in\mathcal{L}^2([0,T]^n)$ is a deterministic symmetric
function; $I_n(h_n)$ stands for the $n$-th order multiple It\^o
integral of $h_n$ with respect to the Brownian motion
$\{B_t\}_{0\leq t\leq T}$. Moreover one has
\begin{eqnarray*}
\Vert X\Vert_2^2&:=&E[|X|^2]\\
&=&\sum_{n\geq 0}n!|h_n|^2_{\mathcal{L}^2([0,T]^n)}.
\end{eqnarray*}

We now introduce the following family of Hilbert spaces of smooth
random variables that was defined by Potthoff and Timpel \cite{PT}
and further studied by Benth and Potthoff \cite{BP}.\\
For $\lambda\geq 1$, let
\begin{eqnarray*}
\mathcal{G}_{\lambda}&:=&\Big\{X=\sum_{n\geq
0}I_n(h_n)\in\mathcal{L}^2(\Omega): \sum_{n\geq
0}n!\lambda^{2n}|h_n|^2_{\mathcal{L}^2([0,T]^n)}< +\infty\Big\}.
\end{eqnarray*}
Note that $\mathcal{G}_1=\mathcal{L}^2(\Omega)$ and for
$1\leq\lambda<\mu$,
$\mathcal{G}_{\mu}\subset\mathcal{G}_{\lambda}\subset\mathcal{L}^2(\Omega)$.
We also denote
\begin{eqnarray*}
\mathcal{G}& := & \bigcap_{\lambda\geq 1}\mathcal{G}_{\lambda}.
\end{eqnarray*}
The most representative element of $\mathcal{G}$ is $\mathcal{E}(f)$
with $f\in\mathcal{L}^2([0,T])$. In fact, since
\begin{eqnarray*}
\mathcal{E}(f)=\sum_{n\geq 0}I_n\Big(\frac{f^{\otimes n}}{n!}\Big),
\end{eqnarray*}
for any $\lambda\geq 1$ one has
\begin{eqnarray*}
\sum_{n\geq
0}\frac{\lambda^{2n}}{n!}|f|^{2n}_{\mathcal{L}^2([0,T])}< +\infty.
\end{eqnarray*}

If $A:\mathcal{L}^2([0,T]) \to \mathcal{L}^2([0,T])$ is a bounded
linear operator then its \emph{second quantization operator}
$\Gamma(A):\mathcal{G} \to \mathcal{G}$ is defined as
\begin{eqnarray*}
\Gamma(A)X & = & \Gamma(A)\sum_{n\geq 0}I_n(h_n)\\
& := & \sum_{n\geq 0}I_n(A^{\otimes n}h_n).
\end{eqnarray*}
With this notation the space $\mathcal{G}_{\lambda}$ previously
defined can be described as
\begin{eqnarray*}
\mathcal{G}_{\lambda}=\{X\in\mathcal{L}^2(\Omega): \Vert
X\Vert_{\mathcal{G}_{\lambda}}:=\Vert\Gamma(\lambda
I)X\Vert_2<+\infty\},
\end{eqnarray*}
where $I$ stands for the identity operator on
$\mathcal{L}^2([0,T])$. In the sequel the operator $\Gamma(\lambda
I)$ will be denoted
simply by $\Gamma(\lambda)$.\\
Observe that
\begin{eqnarray*}
\Gamma(A)\mathcal{E}(f)=\mathcal{E}(Af).
\end{eqnarray*}

Finally let $X=\sum_{n\geq 0}I_n(h_n)\in\mathcal{L}^2(\Omega)$. For
$t\in [0,T]$ the random variable:
\begin{eqnarray*}
D_tX:=\sum_{n\geq 1}nI_{n-1}(h_n(\cdot,t)),
\end{eqnarray*}
where $h_n(\cdot,t)$ is now considered as a function of $n-1$
variables, is called the \emph{Hida-Malliavin derivative} of $X$ at
$t$. By iteration we also define for $k\geq 2$ the $k$-th order
Hida-Malliavin derivative of $X$ at $(t_1,...,t_k)\in[0,T]^k$ as
\begin{eqnarray*}
D^k_{t_1,...,t_k}X:=\sum_{n\geq
k}n(n-1)\cdot\cdot\cdot(n-k+1)I_{n-k}(h_n(\cdot,t_1,...,t_k)).
\end{eqnarray*}
It is easy to see that if $X\in\mathcal{G}$ then for any $k\geq 1$
and any $(t_1,...,t_k)\in[0,T]^k$ the random variable
$D^k_{t_1,...,t_k}X$ belongs to $\mathcal{L}^2(\Omega)$ and
\begin{eqnarray*}
E\Big[\int_{[0,T]^k}|D^k_{t_1,...,t_k}X|^2dt_1...dt_k\Big]<+\infty.
\end{eqnarray*}
A direct calculation shows that for $k\geq 1$,
\begin{eqnarray*}
D^k_{t_1,...,t_k}\mathcal{E}(f)=f(t_1)\cdot\cdot\cdot
f(t_k)\mathcal{E}(f).
\end{eqnarray*}

\section{ A new product for smooth random variables }
We begin introducing a family of products.

\begin{definition}
Let $\varphi:\mathbb{R}\to\mathbb{R}$ be a real analytic function
such that $\varphi(0)=1$ and with expansion,
\begin{eqnarray*}
\varphi(x)=\sum_{k\geq0}a_k x^k,\quad x\in\mathbb{R}.
\end{eqnarray*}
For $X, Y \in \mathcal{G}$ define their
\emph{$\circ_{\varphi}$-product} as
\begin{eqnarray*}
X\circ_\varphi
Y:=\sum_{n\geq0}a_n\int_{[0,T]^n}D_{t_1,...,t_n}^nX\cdot
D_{t_1,...,t_n}^nYdt_1...dt_n.
\end{eqnarray*}
\end{definition}

To ease the notation we will write from now on
\begin{eqnarray*}
\int_{[0,T]^n}D_t^nX\cdot D^n_tYdt
\end{eqnarray*}
to denote the quantity
\begin{eqnarray*}
\int_{[0,T]^n}D_{t_1,...,t_n}^nX\cdot D_{t_1,...,t_n}^nYdt_1...dt_n.
\end{eqnarray*}

\begin{remark}
Observe that for $\varphi(x)=1$ we get $X\circ_{\varphi}Y=X\cdot Y$
while for $\varphi(x)=e^{-x}$ we get $X\circ_\varphi Y=X\diamond Y$
(for a discussion of the last identity see \cite{HY},\cite{L} and
the references quoted there).
\end{remark}

The $\circ_{\varphi}$-product is clearly commutative and
distributive with respect to the sum. We have required the condition
$\varphi(0)=1$ so that $X\circ_\varphi Y$ reduces to $X\cdot Y$ in
the case where $X$ or $Y$ is constant.\\
The next theorem provides a necessary condition on $\varphi$ for the
associativity of the corresponding $\circ_{\varphi}$-product.

\begin{theorem}
If the $\circ_\varphi$-product is associative then $\varphi(x)=
e^{\alpha x}$ for some $\alpha\in\mathbb{R}$.
\end{theorem}

\begin{proof}
Assume that the $\circ_\varphi$-product is associative; this means
that for any $X,Y,Z\in\mathcal{G}$ such that $X\circ_\varphi Y,
Y\circ_\varphi Z\in\mathcal{G}$, the equality
\begin{equation}\label{assoc}
(X\circ_\varphi Y)\circ_\varphi Z=X\circ_\varphi (Y\circ_\varphi Z)
\end{equation}
holds true.\\
Fix $f,g,h\in\mathcal{L}^2([0,T])$ and choose
\begin{eqnarray*}
X=\mathcal{E}(f), Y=\mathcal{E}(g)\mbox{ and }Z=\mathcal{E}(h).
\end{eqnarray*}
With these choices the left-hand side of (\ref{assoc}) becomes after
some simple calculations ($\langle\cdot,\cdot\rangle$ denotes the
inner product in $\mathcal{L}^2([0,T])$):
\begin{eqnarray*}
(\mathcal{E}(f)\circ_{\varphi}\mathcal{E}(g))\circ_{\varphi}\mathcal{E}(h)&=&\varphi(\langle
f,g\rangle)(\mathcal{E}(f)\mathcal{E}(g))\circ_{\varphi}\mathcal{E}(h)\\
&=&\varphi(\langle
f,g\rangle)(\mathcal{E}(f+g)e^{\langle f,g\rangle})\circ_{\varphi}\mathcal{E}(h)\\
&=&\varphi(\langle f,g\rangle)\varphi(\langle
f+g,h\rangle)e^{\langle f,g\rangle}\mathcal{E}(f+g)\mathcal{E}(h)\\
&=&\varphi(\langle f,g\rangle)\varphi(\langle
f+g,h\rangle)\mathcal{E}(f)\mathcal{E}(g)\mathcal{E}(h).
\end{eqnarray*}
while proceeding as before the right-hand side reduces to:
\begin{equation*}
\mathcal{E}(f)\circ_\varphi(\mathcal{E}(g)\circ_\varphi
\mathcal{E}(h))=\varphi(\langle f,g+h\rangle)\varphi(\langle
g,h\rangle)\mathcal{E}(f)\mathcal{E}(g)\mathcal{E}(h).
\end{equation*}
Therefore equation (\ref{assoc}) now reads
\begin{equation*}
\varphi(\langle f,g\rangle)\varphi(\langle
f+g,h\rangle)=\varphi(\langle f,g+h\rangle)\varphi(\langle
g,h\rangle),
\end{equation*}
or equivalently
\begin{equation*}
\varphi(\langle f,g\rangle)\varphi(\langle f,h\rangle+\langle
g,h\rangle)=\varphi(\langle f,g\rangle+\langle
f,h\rangle)\varphi(\langle g,h\rangle),
\end{equation*}
Since $f,g,h\in\mathcal{L}^2([0,T])$ were fixed but arbitrary we can
choose $f$ and $g$ to be orthogonal in $\mathcal{L}^2([0,T])$ and
obtain the equation (recall that by assumption $\varphi(0)=1$),
\begin{eqnarray*}
\varphi(\langle f,h\rangle+\langle g,h\rangle)=\varphi(\langle
f,h\rangle)\varphi(\langle g,h\rangle).
\end{eqnarray*}
Since the equation above is satisfied only by exponential functions,
among the class of continuous functions, the proof is complete.
\end{proof}

We now focus our attention on one specific $\circ_{\varphi}$-product
which will be shown to be related to the heat equation. We begin by
proving a regularity result.

\begin{lemma}
Let $X,Y\in\mathcal{G}_{\sqrt{2}}$. Then the series
\begin{equation*}
\sum_{n\geq0}\frac{1}{n!}\int_{[0,T]^n} D_t^nX\cdot D_t^nY dt,
\end{equation*}
converges in $\mathcal{L}^1(\Omega)$. More precisely,
\begin{eqnarray}\label{estimate}
\Big\Vert\sum_{n\geq0}\frac{1}{n!}\int_{[0,T]^n} D_t^nX\cdot D_t^nY
dt\Big\Vert_1\leq \Vert X\Vert_{\mathcal{G}_{\sqrt{2}}}\Vert
Y\Vert_{\mathcal{G}_{\sqrt{2}}}.
\end{eqnarray}
\end{lemma}

\begin{proof}
Denote by $\Vert\cdot\Vert_1$ the norm in $\mathcal{L}^1(\Omega)$.
By means of the triangle and Cauchy-Schwarz inequalities we get
\begin{eqnarray*}
\Big\Vert\sum_{n\geq0}\frac{1}{n!}\int_{[0,T]^n} D_t^nX\cdot D_t^nY
dt\Big\Vert_1&\leq&\sum_{n\geq0}\frac{1}{n!}\Big\Vert\int_{[0,T]^n}
D_t^nX\cdot
D_t^nYdt\Big\Vert_1\\
&\leq&\sum_{n\geq0}\frac{1}{n!}\int_{[0,T]^n}
\Vert D_t^nX\cdot D_t^nY\Vert_1dt\\
&\leq&\sum_{n\geq0}\frac{1}{n!}\int_{[0,T]^n}
\Vert D_t^nX\Vert_2\cdot\Vert D_t^nY\Vert_2dt\\
&\leq&\sum_{n\geq0}\frac{1}{n!}\Big(\int_{[0,T]^n} \Vert
D_t^nX\Vert^2_2dt\Big)^{\frac{1}{2}}\\
&&\times\Big(\int_{[0,T]^n}
\Vert D_t^nY\Vert^2_2dt\Big)^{\frac{1}{2}}\\
&\leq&\Big(\sum_{n\geq0}\frac{1}{n!}\int_{[0,T]^n} \Vert
D_t^nX\Vert^2_2dt\Big)^{\frac{1}{2}}\\
&&\times\Big(\sum_{n\geq0}\frac{1}{n!}\int_{[0,T]^n} \Vert
D_t^nY\Vert^2_2dt\Big)^{\frac{1}{2}}.
\end{eqnarray*}
Let us now consider the quantity
\begin{eqnarray*}
\sum_{n\geq0}\frac{1}{n!}\int_{[0,T]^n} \Vert D_t^nX\Vert^2_2dt.
\end{eqnarray*}
If $\sum_{k\geq 0}I_k(h_k)$ is the Wiener-It\^o chaos decomposition
of $X$, then
\begin{eqnarray*}
\int_{[0,T]^n} \Vert D_t^nX\Vert^2_2dt=\sum_{k\geq
n}\frac{k!^2}{(k-n)!}|h_k|^2_{\mathcal{L}^2([0,T]^k)},
\end{eqnarray*}
and hence
\begin{eqnarray*}
\sum_{n\geq0}\frac{1}{n!}\int_{[0,T]^n} \Vert
D_t^nX\Vert^2_2dt&=&\sum_{n\geq0}\frac{1}{n!}\sum_{k\geq
n}\frac{k!^2}{(k-n)!}|h_k|^2_{\mathcal{L}^2([0,T]^k)}\\
&=&\sum_{k\geq
0}k!|h_k|^2_{\mathcal{L}^2([0,T]^k)}\sum_{n=0}^k\frac{k!}{n!(k-n)!}\\
&=&\sum_{k\geq
0}k!2^k|h_k|^2_{\mathcal{L}^2([0,T]^k)}\\
&=&E[|\Gamma(\sqrt{2})X|^2]\\
&=&\Vert X\Vert^2_{\mathcal{G}_{\sqrt{2}}}.
\end{eqnarray*}
The same reasoning can be carried for $Y$ completing the proof.
\end{proof}

In view of Remark 3.2 and Lemma 3.4 we make the following
definition.

\begin{definition}
Let $X,Y\in\mathcal{G}_{\sqrt{2}}$. The \emph{anti-Wick product} of
$X$ and $Y$, denoted by $X\circ Y$, is the element of
$\mathcal{L}^1(\Omega)$ defined as
\begin{equation*}
X\circ Y:=\sum_{n\geq0}\frac{1}{n!}\int_{[0,T]^n} D_t^nX\cdot D_t^nY
dt.
\end{equation*}
\end{definition}

\begin{remark}
The anti-Wick product corresponds to the $\circ_\varphi$-product
with $\varphi(x)=e^x$.
\end{remark}

We now prove a crucial result.

\begin{theorem}
Let $X,Y\in\mathcal{G}_{\sqrt{2}}$. Then we have the representation,
\begin{eqnarray*}
X\circ
Y=\Gamma\Big(\frac{1}{\sqrt{2}}\Big)(\Gamma(\sqrt{2})X\cdot\Gamma(\sqrt{2})Y).
\end{eqnarray*}
\end{theorem}

\begin{proof}
First of all we observe that due to inequality (\ref{estimate}), if
$X_n$ converges to $X$ in $\mathcal{G}_{\sqrt{2}}$, then $X_n\circ
Y$ converges to $X\circ Y$ in $\mathcal{L}^1(\Omega)$. Moreover the
family of stochastic exponentials
\begin{eqnarray*}
\mathcal{E}(f)=\exp\Big\{\int_0^Tf(s)dB_s-\frac{1}{2}\int_0^Tf^2(s)ds\Big\},
f\in\mathcal{L}^2([0,T]),
\end{eqnarray*}
forms a total set in $\mathcal{G}_{\sqrt{2}}$. Therefore by
linearity and continuity it is sufficient to prove the theorem for
stochastic exponentials. But this is easily done; indeed for any
$f,g\in\mathcal{L}^2([0,T])$,
\begin{eqnarray*}
\mathcal{E}(f)\circ\mathcal{E}(g)&=&\sum_{n\geq0}\frac{1}{n!}\int_{[0,T]^n}
D_t^n\mathcal{E}(f)\cdot D_t^n\mathcal{E}(g)dt\\
&=&\mathcal{E}(f)\mathcal{E}(g)e^{\int_0^Tf(s)g(s)ds},
\end{eqnarray*}
and
\begin{eqnarray*}
\Gamma\Big(\frac{1}{\sqrt{2}}\Big)(\Gamma(\sqrt{2})\mathcal{E}(f)\cdot\Gamma(\sqrt{2})\mathcal{E}(g))&=&
\Gamma\Big(\frac{1}{\sqrt{2}}\Big)(\mathcal{E}(\sqrt{2}f)\mathcal{E}(\sqrt{2}g))\\
&=&\Gamma\Big(\frac{1}{\sqrt{2}}\Big)\mathcal{E}(\sqrt{2}(f+g))e^{2\int_0^Tf(s)g(s)ds}\\
&=&\mathcal{E}(f+g)e^{2\int_0^Tf(s)g(s)ds}\\
&=&\mathcal{E}(f)\mathcal{E}(g)e^{\int_0^Tf(s)g(s)ds}.
\end{eqnarray*}
\end{proof}

\begin{corollary}
The anti-Wick product is associative.
\end{corollary}
\begin{proof}
Let $X,Y,Z\in\mathcal{G}_{\sqrt{2}}$ be such that $X\circ Y, Y\circ
Z\in\mathcal{G}_{\sqrt{2}}$. We have to prove that
\begin{eqnarray*}
(X\circ Y)\circ Z=X\circ (Y\circ Z).
\end{eqnarray*}
By the representation of Theorem 3.7, this follows immediately by a
straightforward verification.
\end{proof}

\section{ Application to the heat equation }

The representation in Theorem 3.7 enables us to write for
$X\in\mathcal{G}_{\sqrt{2}}$ and $n\geq 1$,
\begin{eqnarray*}
X^{\circ n}&:=&\underbrace{X\circ\cdot\cdot\cdot\circ X}_{n-times}\\
&=&\Gamma\Big(\frac{1}{\sqrt{2}}\Big)\Big((\Gamma(\sqrt{2})X)^n\Big).
\end{eqnarray*}
Therefore if $f:\mathbb{R}\to\mathbb{R}$ is such that
$f(\Gamma(\sqrt{2})X)\in\mathcal{L}^1(\Omega)$ we can define
\begin{eqnarray*}
f^{\circ}(X):=\Gamma\Big(\frac{1}{\sqrt{2}}\Big)f(\Gamma(\sqrt{2})X),
\end{eqnarray*}
as an element of $\mathcal{L}^1(\Omega)$.
\begin{example}
Let $f(x)=\exp\{x\}$ and $h\in\mathcal{L}^2([0,T])$. Then
\begin{eqnarray*}
\exp^{\circ}\Big\{\int_0^Th(s)dB_s\Big\}&=&\Gamma\Big(\frac{1}{\sqrt{2}}\Big)\exp\Big
\{\Gamma(\sqrt{2})\int_0^Th(s)dB_s\Big\}\\
&=&\Gamma\Big(\frac{1}{\sqrt{2}}\Big)\exp\Big
\{\sqrt{2}\int_0^Th(s)dB_s\Big\}\\
&=&\Gamma\Big(\frac{1}{\sqrt{2}}\Big)\mathcal{E}(\sqrt{2}h)e^{\int_0^Th^2(s)ds}\\
&=&\mathcal{E}(h)e^{\int_0^Th^2(s)ds}\\
&=&\exp\Big\{\int_0^Th(s)dB_s+\frac{1}{2}\int_0^Th(s)ds\Big\}.
\end{eqnarray*}
Observe the symmetry with the case
\begin{eqnarray*}
\exp^{\diamond}\Big\{\int_0^Th(s)dB_s\Big\}=\exp\Big\{\int_0^Th(s)dB_s-\frac{1}{2}\int_0^Th(s)ds\Big\}.
\end{eqnarray*}

\end{example}
We now come to the main theorems of the present section which
establish a new probabilistic representation for the solution of the
heat equation in terms of anti-Wick products.

\begin{theorem}
Let $f:\mathbb{R}\to\mathbb{R}$ be a bounded and continuous function
and let $u:[0,T]\times\mathbb{R}\to\mathbb{R}$ be the unique
solution of
\begin{eqnarray*}
\Big{ \{ }
\begin{array}{ll}
\partial_tu(t,x)=\frac{1}{2}\partial_{xx} u(t,x),\quad x\in\mathbb{R},t\in]0,T]\\
u(0,x)=f(x),\quad x\in\mathbb{R}
\end{array}
\end{eqnarray*}
among the class of functions satisfying the bound
\begin{eqnarray*}
|u(t,x)|\leq Ae^{ax^2},\quad x\in\mathbb{R}, t\in[0,T],
\end{eqnarray*}
for some $a,A>0$. Then for $t\in[0,T]$,
\begin{eqnarray*}
u(t,B_t)=f^{\circ}(B_t).
\end{eqnarray*}
\end{theorem}

\begin{proof}
Recall that
\begin{eqnarray*}
f^{\circ}(B_t)&=&\Gamma\Big(\frac{1}{\sqrt{2}}\Big)f(\Gamma(\sqrt{2})B_t)\\
&=&\Gamma\Big(\frac{1}{\sqrt{2}}\Big)f(\sqrt{2}B_t).
\end{eqnarray*}
Fix $h\in\mathcal{L}^2([0,T])$; we will prove the theorem by showing
that
\begin{eqnarray*}
E[u(t,B_t)\mathcal{E}(h)]=E[f^{\circ}(B_t)\mathcal{E}(h)],
\end{eqnarray*}
for any $h\in\mathcal{L}^2([0,T])$. Using the properties of second
quantization operators and the Girsanov theorem we get
\begin{eqnarray*}
E[f^{\circ}(B_t)\mathcal{E}(h)]&=&E\Big[\Big(\Gamma\Big(\frac{1}{\sqrt{2}}\Big)f(\sqrt{2}B_t)\Big)\mathcal{E}(h)\Big]\\
&=&E\Big[f(\sqrt{2}B_t)\mathcal{E}\Big(\frac{h}{\sqrt{2}}\Big)\Big]\\
&=&E\Big[f\Big(\sqrt{2}B_t+\int_0^th(s)ds\Big)\Big].
\end{eqnarray*}
The law of $\sqrt{2}B_t$ is the same as the one of $B_t+\tilde{B}_t$
where $\tilde{B}_t$ is another Brownian motion independent of $B_t$
and defined on an auxiliary probability space
$(\tilde{\Omega},\tilde{\mathcal{F}},\tilde{\mathcal{P}})$.
Therefore, denoting by $\tilde{E}$ the expectation in this new
probability space, we obtain
\begin{eqnarray*}
E\Big[f\Big(\sqrt{2}B_t+\int_0^th(s)ds\Big)\Big]&=&E\Big[\tilde{E}\Big[f\Big(B_t+\tilde{B}_t+\int_0^th(s)ds\Big)\Big]\Big]\\
&=&E\Big[u\Big(t,B_t+\int_0^th(s)ds\Big)\Big]\\
&=&E[u(t,B_t)\mathcal{E}(h)],
\end{eqnarray*}
where we used the well known formula
\begin{eqnarray*}
u(t,x)=\tilde{E}[f(\tilde{B}_t+x)],\quad x\in\mathbb{R}, t\in [0,T].
\end{eqnarray*}
The proof is complete.
\end{proof}

Under the same assumptions on the initial condition we can also
prove the following.

\begin{theorem}
Let $u(t,x)$ be the unique solution of
\begin{eqnarray*}
\Big{ \{ }
\begin{array}{ll}
\partial_tu(t,x)=\frac{1}{2}\partial_{xx} u(t,x),\quad x\in\mathbb{R},t\in ]0,T],\\
u(0,x)=f(x),\quad x\in\mathbb{R},
\end{array}
\end{eqnarray*}
and $v(t,x)$ be the unique solution of
\begin{eqnarray*}
\Big{ \{ }
\begin{array}{ll}
\partial_tv(t,x)=\frac{1}{2}\partial_{xx} v(t,x),\quad x\in\mathbb{R}, t\in ]0,T]\\
v(0,x)=g(x),\quad x\in\mathbb{R}.
\end{array}
\end{eqnarray*}
Moreover let $w(t,x)$ be the unique solution of
\begin{eqnarray*}
\Big{ \{ }
\begin{array}{ll}
\partial_tw(t,x)=\frac{1}{2}\partial_{xx} w(t,x),\quad x\in\mathbb{R}, t\in ]0,T]\\
w(0,x)=f(x)g(x),\quad x\in\mathbb{R}.
\end{array}
\end{eqnarray*}
Then for each $t\in [0,T]$,
\begin{eqnarray*}
u(t,B_t)\circ v(t,B_t)=w(t,B_t).
\end{eqnarray*}
\end{theorem}

\begin{proof}
By Theorem 3.7 we have that
\begin{eqnarray*}
u(t,B_t)\circ
v(t,B_t)&=&\Gamma\Big(\frac{1}{\sqrt{2}}\Big)(\Gamma(\sqrt{2})u(t,B_t)\cdot\Gamma(\sqrt{2})v(t,B_t))\\
&=&\Gamma\Big(\frac{1}{\sqrt{2}}\Big)(f(\sqrt{2}B_t)g(\sqrt{2}B_t))\\
&=&\Gamma\Big(\frac{1}{\sqrt{2}}\Big)((f\cdot g)(\sqrt{2}B_t))\\
&=&w(t,B_t).
\end{eqnarray*}
\end{proof}

\section{Appendix}
In this section we prove two formulas relating Wick and anti-Wick
products.

\begin{proposition}
Let $X, Y\in\mathcal{G}$. Then,
\begin{eqnarray*}
X\circ Y=\sum_{n\geq0}\frac{2^n}{n!}\int_{[0,T]^n}D_t^nX\diamond
D_t^nY dt.
\end{eqnarray*}
\end{proposition}

\begin{proof}
We have
\begin{eqnarray*}
X\circ Y &=&\sum_{j\geq0}\frac{1}{j!}\int_{[0,T]^j}D_t^j X\cdot D_t^jY dt\\
&=&\sum_{j\geq0}\frac{1}{j!}\int_{[0,T]^j}
\Big(\sum_{k\geq0}\frac{1}{k!}\int_{[0,T]^k}D_{s,t}^{j+k}X
\diamond D_{s,t}^{j+k}Y ds\Big)dt\\
&=&\sum_{j\geq0}\sum_{k\geq
0}\frac{1}{j!}\frac{1}{k!}\int_{[0,T]^j}\int_{[0,T]^k}D_{s,t}^{j+k}X
\diamond D_{s,t}^{j+k}Y ds dt\\
&=&\sum_{n\geq0}\frac{1}{n!}\int_{[0,T]^n}D_t^n X\diamond D_t^n Ydt\sum_{k=0}^n{n \choose k}\\
&=&\sum_{n\geq0}\frac{2^n}{n!}\int_{[0,T]^n}D_t^n X\diamond D_t^nY
dt.
\end{eqnarray*}
\end{proof}

\begin{proposition}
Let $X, Y\in\mathcal{G}$. Then,
\begin{eqnarray*}
X\diamond Y=\sum_{n\geq0}\frac{(-2)^n}{n!}\int_{[0,T]^n}D_t^n X\circ
D_t^nY dt.
\end{eqnarray*}
\end{proposition}
\begin{proof}
For $S,T\in\mathcal{G}$ we can write
\begin{eqnarray*}
S\diamond T&=&\sum_{n\geq0}\frac{(-1)^n}{n!}\int_{[0,T]^n}D_t^nX\cdot D_t^nYdt\\
&=&\sum_{n\geq0}\frac{(-1)^n}{n!}
\int_{[0,T]^n}D_t^n\Gamma(\sqrt{2})\Gamma\left(\frac{1}{\sqrt{2}}\right)X
\cdot D_t^n\Gamma(\sqrt{2})\Gamma\left(\frac{1}{\sqrt{2}}\right)Y dt\\
&=&\sum_{n\geq0}\frac{(-1)^n}{n!}
\int_{[0,T]^n}2^{\frac{n}{2}}\Gamma(\sqrt{2})D_t^n\Gamma\left(\frac{1}{\sqrt{2}}\right)X
\cdot 2^{\frac{n}{2}}\Gamma(\sqrt{2})D_t^n\Gamma\left(\frac{1}{\sqrt{2}}\right)Ydt\\
&=&\sum_{n\geq0}\frac{(-2)^n}{n!}
\Gamma(\sqrt{2})\int_{[0,T]^n}D_t^n\Gamma\left(\frac{1}{\sqrt{2}}\right)X
\circ D_t^n\Gamma\left(\frac{1}{\sqrt{2}}\right)Y dt.
\end{eqnarray*}
Applying the operator $\Gamma\left(\frac{1}{\sqrt{2}}\right)$ to the
first and last terms of the previous chain of equalities we obtain
the desired result by letting
$X:=\Gamma\left(\frac{1}{\sqrt{2}}\right)S$ and
$Y:=\Gamma\left(\frac{1}{\sqrt{2}}\right)T$.
\end{proof}

\begin{center}
{\bf Acknowledgments}
\end{center}

The authors are grateful to Aurel I. Stan for stimulating
discussions.

\end{document}